\documentclass[onefignum]{siamart190516}
\pdfoutput=1
\usepackage{multicol}
\usepackage{amsfonts}
\usepackage{xparse}

\usepackage{graphicx}
\usepackage[noend]{algorithmic}
\usepackage{float}
\floatstyle{plain}
\newfloat{multalg}{ht}{multalg}

\graphicspath{{fig/}}

\def\Figure#1#2{
\begin{center}
\includegraphics{#1}
\caption{#2}
\label{fig:#1}
\end{center}
}

\newcommand{\deq}{\doteq}
\newcommand{\la}{\leftarrow}

\def\fft#1{{\tt fft}(#1)}
\def\crfft#1{{\tt crfft}(#1)} 
\def\ifft#1{{\tt ifft}(#1)}
\def\icrfft#1{{\tt rcfft}(#1)}

\def\dft{\mathop{\rm dft}\nolimits}

\def\mult{{\tt mult}}

\def\Forward{{\tt Forward}}
\def\Backward{{\tt Backward}}

\def\convolve{{\tt convolve}}

\def\forward{{\tt forward1}}
\def\forwardTwo{{\tt forward2}}
\def\forwardInner{{\tt forwardInner}}

\def\backward{{\tt backward1}}
\def\backwardTwo{{\tt backward2}}
\def\backwardInner{{\tt backwardInner}}

\def\convolveX{{\tt convolveX}}

\def\X{{\tt X}}
\def\Centered{{\tt C}}
\def\Hermitian{{\tt H}}

\def\forwardXTwo{{\tt forward2X}}
\def\forwardXInner{{\tt forwardInnerX}}

\def\backwardXTwo{{\tt backward2X}}
\def\backwardXInner{{\tt backwardInnerX}}

\def\param{L,m,q}

\def\forwardCentered{{\tt{} forward2C}}
\def\forwardCenteredInner{{\tt forwardInnerC}}

\def\backwardCentered{{\tt backward2C}}
\def\backwardCenteredInner{{\tt backwardInnerC}}

\def\forwardHermitian{{\tt forward2H}}
\def\forwardHermitianInner{{\tt forwardInnerH}}

\def\backwardHermitian{{\tt backward2H}}
\def\backwardHermitianInner{{\tt backwardInnerH}}

\def\Hf{H}
\def\Hc{\tilde{H}}

\def\convdesc#1#2{#1~is a one-dimensional #2
  convolution. There are $A$ inputs (each of length $L$, to be padded
  to at least length $M$) and $B$ outputs. The convolution uses the
  multiplication operator \mult.}

\def\FBdesc#1#2#3#4{is the #4 #1 transform for residue~#2 when #3.}

\def\Fdesc#1#2#3{\FBdesc{forward}{#1}{\hbox{#2}}{#3}}

\def\Bdesc#1#2#3{\FBdesc{backward}{#1}{\hbox{#2}}{#3}}

\def\deltadesc{Here $\delta_t$ is the Kronecker $\delta_{t,0}$.}

\def\ph{{p_2}} 
\def\phm{{\ph m}}
\def\e1{{e_1}}
\def\mod{\mathop{\rm mod}\nolimits}
\newcommand{\ze}[2]{{\zeta_{#1}^{#2}}}

\renewcommand{\vec}[1]{{\boldsymbol{#1}}}
\def\vf{\vec{f}}
\def\vg{\vec{g}}
\def\vh{\vec{h}}
\def\vF{\vec{F}}
\def\vG{\vec{G}}

\NewDocumentCommand{\seq}{m O{#3} m m m}{%
  \left\{#1_{#2}\right\}_{#3=#4}^{#5}
}

\newcommand{\ceil}[1]{{\left\lceil #1 \right\rceil}}
\newcommand{\floor}[1]{{\left\lfloor #1 \right\rfloor}}
\newcommand{\conj}[1]{{\overline{#1}}}

\def\INPUT{\ENSURE}

\def\papertitle{Hybrid Dealiasing of Complex Convolutions}

\headers{Hybrid Dealiasing of Complex Convolutions}
{Noel Murasko and John C. Bowman}

\title{\papertitle}

\author{Noel Murasko$^*$ \and John C. Bowman\thanks{University of Alberta
(\email{murasko@ualberta.ca}, \email{bowman@ualberta.ca})}}

\ifpdf
\hypersetup{
  pdftitle={\papertitle},
  pdfauthor={Noel Murasko and John C. Bowman}
}
\fi

\def\hyperlabel#1{\label{#1}\hypertarget{#1}}
\def\alglink#1{\hyperlink{alg:#1}{\csname #1\endcsname}}
\def\ForwardIs#1{$\Forward\la $\hyperlink{alg:#1}{\csname #1\endcsname}}
\def\BackwardIs#1{$\Backward\la $\hyperlink{alg:#1}{\csname #1\endcsname}}

\begin{document}

\maketitle

\begin{abstract}
Efficient algorithms for computing linear convolutions based on the fast
Fourier transform are developed. A hybrid approach is described that combines
the conventional practice of explicit dealiasing (explicitly padding the input
data with zeros) and implicit dealiasing (mathematically accounting
for these zero values). The new approach generalizes implicit dealiasing to
arbitrary padding ratios and includes explicit dealiasing as a special case.
Unlike existing implementations of implicit dealiasing, hybrid dealiasing
tailors its subtransform sizes to the convolution geometry.
Multidimensional convolutions are implemented with hybrid dealiasing by
decomposing them into lower-dimensional convolutions.
Convolutions of complex-valued and
Hermitian inputs of equal length are illustrated with pseudocode and
implemented in the open-source {\tt FFTW++} library. Hybrid dealiasing
is shown to outperform explicit dealiasing in one, two, and three dimensions.
\end{abstract}



\section{Introduction}
The convolution theorem provides a method for efficiently computing
circular convolutions using Fast Fourier Transforms (FFTs).
However, many applications require linear convolutions; these
can be computed using the convolution theorem by padding the data
with a sufficient number of zeros to avoid polluting the result with
\emph{aliases}: errors that arise from the lack of periodicity in the input.
The standard practice is to pad the input arrays explicitly with zeros before
taking the FFT. While this accomplishes dealiasing, it requires
reading and multiplying values that are known {\it a priori\/} to be zero.

{\it Implicit dealiasing} \cite{Bowman11,Roberts18} provides an
alternative to explicit dealiasing. The FFTs are formulated to take
account of the known zero values implicitly, avoiding the need for
explicit zero padding.
In previous work, the focus was on developing implicit dealiasing for $1/2$
and $2/3$ padding ratios \cite{Orszag71} (these fractions refer to
the ratio of input data length to the zero-padded buffer size) for complex and
Hermitian symmetric input data. While these cases are important for applications,
such as signal processing and pseudospectral methods for solving partial
differential equations \cite{Patterson71,GottliebOrszag77}, many applications
do not satisfy these requirements.

The goal of this work is to formulate and develop a systematic
framework, which we call {\it hybrid dealiasing,} for the efficient dealiasing
of complex-valued convolutions with arbitrary padding ratios.
The decomposition of multidimensional dealiased convolutions into
padded FFTs and one-dimensional subconvolutions \cite{Bowman11,Roberts18}
is crucial for the enhanced performance afforded by hybrid dealiasing,
even when the zero padding is explicit.

In~\cref{sec:standardCase}, we review the concept of implicit
dealiasing and introduce hybrid padding. \Cref{sec:centeredCase,sec:hermitianCase} then extend hybrid padding to centered and
Hermitian symmetric data. In \cref{sec:multidim}, we discuss how to
compute multidimensional convolutions efficiently by decomposing them
into lower dimensional convolutions. An overview of numerical
optimizations and results are given in
\cref{sec:optimization,sec:results}. Finally, future
applications of this work are discussed in \cref{sec:applications}.

\section{Complex convolutions in one dimension}\label{sec:standardCase}
Let~$X$ denote the space of sequences with bounded support equipped
with addition and scalar multiplication.

Let~$A,B\in \mathbb{N}$. An operator~${\cal C}:X^A\mapsto X^B$ is a
\emph{general convolution} if there exists~$M\in \mathbb{N}$ and a
pointwise operator~${\cal M}:X^A\mapsto X^B$
(called the \emph{multiplication operator}), which satisfies
$$
\dft[{\cal C}(\vf_1,\ldots, \vf_A)]={\cal M}(\dft[\vf_1],\ldots,\dft[\vf_A]),
$$
where~$\vf_1,\ldots,\vf_A$ have been zero-padded to length~$M$
and $\dft$ denotes the discrete Fourier transform (DFT) of each
component of its argument.
An important special case is the
convolution of two sequences
$\vf=\{f_{j}\}_{j\in \mathbb{Z}}$ and $\vg=\{g_{j}\}_{j\in \mathbb{Z}}$
given by $(\vf*\vg)_{j}\deq\sum_{i\in\mathbb{Z}}f_{i}g_{j-i}$,
where we use $\deq$ to denote a definition.
One can compute such convolutions efficiently by taking
taking the DFT of each sequence, multiplying the results
element by element, and then taking the inverse DFT.
Due to the cyclic nature of the DFT, the resulting
convolution is circular; for a linear convolution, the input data must
be padded (\emph{dealiased}) with zeros.

The amount of
zero padding needed to compute the convolution using DFTs depends on
the pointwise multiplication operator in the convolution and the
length of the input data. For example, when ${\cal
  M}(\vf_1\ldots,\vf_A)_j=(\vf_1)_j\cdots(\vf_A)_j$, one can obtain
all components of a convolution of $A$ inputs with lengths
$L_1,\ldots,L_A$ by padding each input to size $M=\sum_{i=1}^n L_i - (A-1)$.
However, in some applications not all components of the
convolution are needed; this can reduce the padding requirements.

In this work, we only consider convolutions with $A$ equal-size inputs
and $B$ outputs that result from an arbitrary
element-by-element multiplication in the transformed domain. Given
input data $\{a_j\}_{j=0}^{L-1}$ that needs to be
padded with zeros to length $M$, we construct a buffer $\vf\deq
\{f_j\}_{j=0}^{M-1}$ where $f_j=a_j$ for $j < L$ and $f_j = 0$ for
$j \geq L$. The DFT of $\vf$ can be written as
\begin{equation*}
  F_{k} = \sum_{j=0}^{M-1}\ze{M}{kj}f_{j} = \sum_{j=0}^{L-1}\ze{M}{kj}f_{j},
  \ \ \ k \in \{ 0, \ldots, M -1\},
\end{equation*}
where $\ze{N}{} \deq \exp\left(2\pi i/N\right)$ is the $N$th primitive root of
unity.

For now, assume that $L$ and $M$ share a common factor $m$, so that $L =pm$ and
$M=qm$, where $m, p, q \in \mathbb{N}$, with $q \geq p$. We can then reindex $j$
and $k$ as
$$j = tm+s, \ \ t\in\{0,
\ldots,p-1\}, \ \ s\in\{0,\ldots, m-1\},$$
$$k = q\ell+r, \ \ \ell\in\{0, \ldots, m-1\}, \ \ r\in\{0,\ldots, q-1\}.$$
This allows us to decompose the DFT via the Cooley--Tukey
algorithm \cite{Cooley65} as in \cite{Bowman11}:
\begin{equation}\label{eqn:forStan}
  F_{q\ell+r} = \sum_{s=0}^{m-1}\sum_{t=0}^{p-1}\ze{qm}{(q\ell+r)(tm+s)}f_{tm+s}
  =\sum_{s=0}^{m-1}\ze{m}{\ell s}\ze{qm}{rs}\sum_{t=0}^{p-1}\ze{q}{rt}f_{tm+s}.
\end{equation}
Computing the DFT of $\vf$ then amounts to preprocessing and computing $q$ DFTs of size $m$; no explicit zero
padding is needed. The inverse transform is similar \cite{Bowman11}:
\begin{equation}\label{eqn:invStan}
  f_{tm+s} = \frac{1}{qm}\sum_{r=0}^{q-1}\ze{q}{-tr}\ze{qm}{-sr}
  \sum_{\ell=0}^{m-1}\ze{m}{-s\ell}F_{q\ell+r}.
\end{equation}
This transform requires $q$ DFTs of size $m$, followed by post-processing.
We now demonstrate how these equations can be implemented to
generalize implicit dealiasing to a wider class of convolutions.

\subsection{Hybrid padding}
An issue with the above formulation is the assumption that $L$ and $M$ must
share a common factor $m$. Furthermore, even if $L$ and $M$ do share a common
factor, the resulting convolution might be inefficient: the most efficient
FFTs currently available are for products of powers of the radices
2, 3, 5, and 7 \cite{Frigo05}.

\begin{figure}[htbp]
\begin{center}
\includegraphics[width=\linewidth]{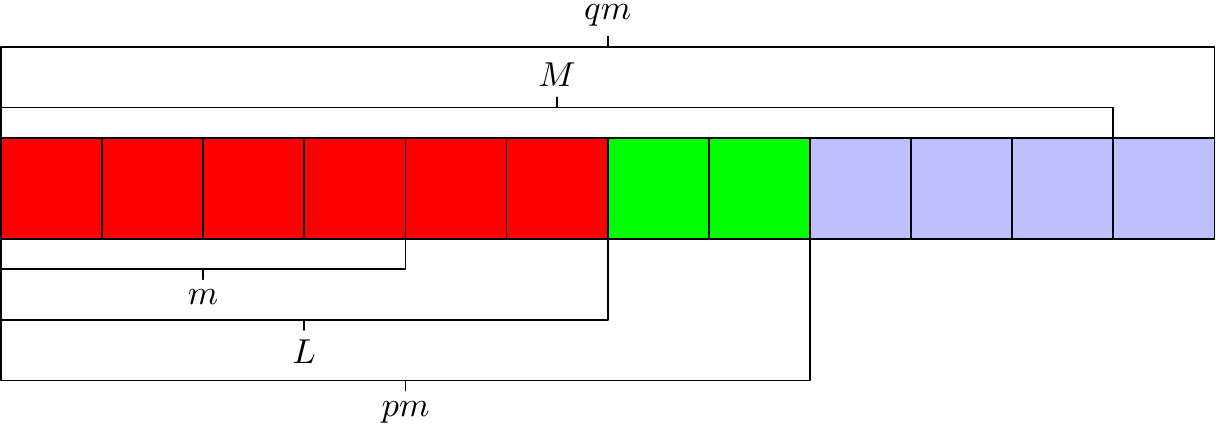}
\end{center}
\caption{An illustration of hybrid padding for a one-dimensional
  array with $L=6$ and $M=11$. Choosing $m=4$, we have $p=2$ and
  $q=3$. We explicitly pad our data of length $L$ to length $pm=8$,
  and then implicitly pad our data to length $qm=12$.}
\label{fig:hybridPadding}
\end{figure}

Our solution to this problem relies on the observation that
$M$ is the \emph{minimum} size required to dealias a convolution:
padding beyond $M$ is fine (and perhaps desired if it increases
efficiency). Given some $m \in \mathbb{N}$, we define
\begin{equation}\label{eqn:pqdef1}
  p\deq\ceil{\frac{L}{m}}, \ \ \ \ q\deq\ceil{\frac{M}{m}}.
\end{equation}
These are the smallest positive integers such that $pm \geq L$ and $qm \geq
M$. To take the forward transform, we \emph{explicitly} pad $\vf$ with zeros to
size $pm$ and then use \cref{eqn:forStan} to compute the padded transform of
size $qm$. To take the inverse transform, we use \cref{eqn:invStan},
ignoring the last $pm-L$ elements. We refer to this combination of explicit and
implicit padding, illustrated in \cref{fig:hybridPadding}, as
\emph{hybrid padding}. If out-of-place FFTs are used, any explicit zero
padding only needs to be written to the buffer once.

An advantage of hybrid padding is the ability to choose any $m$ value, as the
choice of $m$ is independent of the size of the input array.
For the remainder of this work, we exclusively refer to padding from size $pm$
to $qm$, keeping in mind that we might have to use hybrid padding to achieve
this.

\subsection{Convolutions one residue at a time}\label{sub:oneRes}

\begin{multalg}[htbp]
\begin{minipage}[t]{0.49\textwidth}
\centering
\begin{algorithm}[H]
\caption{\convdesc{\convolve}{complex}}
\begin{algorithmic}\hyperlabel{alg:convolve}
\INPUT{$\seq{\vf}{a}{0}{A-1}$, $L$, $M$, $m$, $A$, $B$}
\STATE{$p \la \ceil{L/m}$}
\IF{$p=1$}
  \STATE{$n\la \ceil{M/m}$}
  \STATE{$q\la n$}
  \STATE{\ForwardIs{forward}}
  \STATE{\BackwardIs{backward}}
\ELSIF{$p=2$}
  \STATE{$n\la \ceil{M/m}$}
  \STATE{$q\la n$}
  \STATE{\ForwardIs{forwardTwo}}
  \STATE{\BackwardIs{backwardTwo}}
\ELSE
  \STATE{$n\la \ceil{M/m}$}
  \STATE{$q\la np$}
  \STATE{\ForwardIs{forwardInner}}
  \STATE{\BackwardIs{backwardInner}}
\ENDIF
\FOR{$b = 0, \ldots, B-1$}
  \STATE{$\vh_b \la \seq{0}[]{j}{0}{L-1}$}
\ENDFOR
\FOR{$r = 0, \ldots, n-1$}
  \FOR{$a = 0, \ldots, A-1$}
    \STATE{$\vF_a\la \Forward(\vf_a,\param,r)$}
  \ENDFOR
  \STATE{$\seq{\vF}{b}{0}{B-1}\la \mult(\seq{\vF}{a}{0}{A-1})$}
  \FOR{$b = 0, \ldots, B-1$}
    \STATE{$\vh_b \la \vh_b+\Backward(\vF_b,\param,r)$}
  \ENDFOR
\ENDFOR
\FOR{$b = 0, \ldots, B-1$}
  \STATE{$\vf_b\la \vh_b/(qm)$}
\ENDFOR
\RETURN{$\seq{\vf}{b}{0}{B-1}$}
\end{algorithmic}
\end{algorithm}
\end{minipage}
\hfill
\begin{minipage}[t]{0.49\textwidth}
\centering
\begin{algorithm}[H]
\caption{\convdesc{\convolveX}{centered or Hermitian}~In the centered
version, \X~denotes \Centered. In the Hermitian version, \X~denotes
\Hermitian~and only $\ceil{L/2}$ inputs are provided (corresponding to the
non-negative indices).}
\begin{algorithmic}\hyperlabel{alg:convolveCentered}
\INPUT{$\seq{\vf}{a}{0}{A-1}$, $L$, $M$, $m$, $A$, $B$}
\STATE{$p \la 2\ceil{L/(2m)}$}\ \COMMENT{force even $p$}
\IF{$p=2$}
  \STATE{$\Forward\la $\forwardXTwo}
  \STATE{$\Backward\la $\backwardXTwo}
\ELSE
  \STATE{$\Forward\la $\forwardXInner}
  \STATE{$\Backward\la $\backwardXInner}
\ENDIF
\STATE{$n\la \ceil{2M/(pm)}$}
\STATE{$q\la np/2$}
\FOR{$b = 0, \ldots, B-1$}
  \STATE{$\vh_b \la \seq{0}[]{j}{0}{L-1}$}
\ENDFOR
\FOR{$r = 0, \ldots, n-1$}
  \FOR{$a = 0, \ldots, A-1$}
    \STATE{$\vF_a\la\Forward(\vf_a,\param, r)$}
  \ENDFOR
  \STATE{$\seq{\vF}{b}{0}{B-1}\la \mult(\seq{\vF}{a}{0}{A-1})$}
  \FOR{$b = 0, \ldots, B-1$}
    \STATE{$\vh_b\la \vh_b+\Backward(\vF_b,\param, r)$}
  \ENDFOR
\ENDFOR
\FOR{$b = 0, \ldots, B-1$}
  \STATE{$\vf_b\la \vh_b/(qm)$}
\ENDFOR
\RETURN{$\seq{\vf}{b}{0}{B-1}$}
\end{algorithmic}
\end{algorithm}
\end{minipage}
\end{multalg}

\begin{figure}[htbp]
\begin{center}
\includegraphics[width=\linewidth]{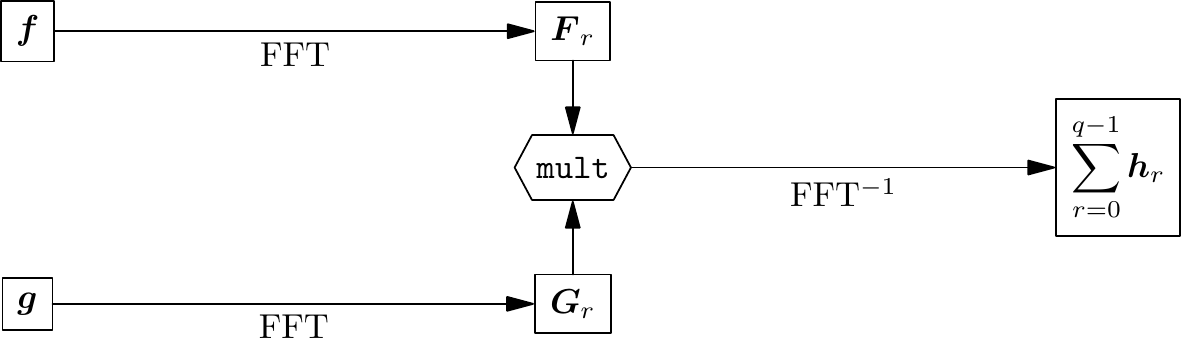}
\caption{Accumulation of residue contributions to a convolution.}
\label{fig:accumulation}
\end{center}
\end{figure}

For each $r\in\{0,\ldots q-1\}$, we define the \emph{residue contribution}
$\vF_r\doteq\seq{F}[q\ell + r]{\ell}{0}{m-1}$. A key optimization that
allows us to save memory is that we can
compute contributions to the convolution one residue at a time. To find the
inverse for that residue contribution, we compute
\begin{equation*}
h_{r,s,t}\doteq \ze{q}{-tr}\ze{qm}{-sr}\sum_{\ell=0}^{m-1}
\ze{m}{-s\ell}F_{q\ell+r}.
\end{equation*}
Accumulating over $r\in\{0,\ldots, q-1\}$, we obtain the inverse:
\begin{equation*}
f_{tm+s} = \frac{1}{qm}\sum_{r=0}^{q-1}h_{r,s,t}.
\end{equation*}
This formulation can be advantageous for large problems as it allows
reuse of the memory needed to store $\vF_r$, illustrated
for a binary convolution in \cref{fig:accumulation}.

A general convolution applies a given multiplication operator to
the DFT of $A$ inputs to produce $B$ outputs, which are then
inverse transformed to produce the final result. The contributions to the
forward transform from each residue can either be computed independently or
grouped together. When using the conjugate symmetry optimization (described in
\cref{sub:conjSymOpt}), it is convenient to compute the~$r$ and~$-r\mod q$
residues together, with the $r=0$ and $r=q/2$ residues
(if~$q$ is even) stored together.

Pseudocode for an in-place convolution is given in
\cref{alg:convolve}.
Pseudocode for the forward and backward
transforms is shown in \cref{alg:forward,alg:backward}
for $p=1$ and \cref{alg:forwardTwo,alg:backwardTwo}
for $p=2$. For $p > 2$, it is more efficient to replace the
summations in \cref{eqn:forStan,eqn:invStan}
by an inner DFT, which is discussed in \cref{sub:inner}.
The pseudocode illustrates only the case where one
residue is computed at a time.
\begin{multalg}[htbp]
\begin{minipage}[t]{0.49\textwidth}
\begin{algorithm}[H]
\caption{\forward~\Fdesc{$r$}{$p=1$}{complex}}
\begin{algorithmic}\hyperlabel{alg:forward}
\INPUT{$\seq{f}{j}{0}{L-1},\param,r$}
\FOR{$s = 0, \ldots, L-1$}
  \STATE{$W_s \la \ze{qm}{rs}f_s$ }
\ENDFOR
\FOR{$s = L, \ldots, m-1$}
  \STATE{$W_s \la 0$ }
\ENDFOR
\STATE{$\seq{V}{\ell}{0}{m-1}\la \fft{\seq{W}{s}{0}{m-1}}$}
\RETURN{$\seq{V}{k}{0}{m-1}$}
\end{algorithmic}
\end{algorithm}
\end{minipage}
\hfill
\begin{minipage}[t]{0.49\textwidth}
\begin{algorithm}[H]
\caption{\backward~\Bdesc{$r$}{$p=1$}{complex}}
\begin{algorithmic}\hyperlabel{alg:backward}
\INPUT{$\seq{F}{k}{0}{m-1},\param,r$}
\STATE{$\seq{W}{s}{0}{m-1}\la \ifft{\seq{F}{\ell}{0}{m-1}}$}
\FOR{$s = 0, \ldots, L-1$}
  \STATE{$W_s \la \ze{qm}{-rs}W_s$}
\ENDFOR
\RETURN{$\seq{W}{j}{0}{L-1}$}
\end{algorithmic}
\end{algorithm}
\end{minipage}
\end{multalg}

\begin{multalg}[htbp]
\begin{minipage}[t]{0.49\textwidth}
\centering
\begin{algorithm}[H]
\caption{\forwardTwo~\Fdesc{$r$}{$p=2$}{complex}}
\begin{algorithmic}\hyperlabel{alg:forwardTwo}
\INPUT{$\seq{f}{j}{0}{L-1},\param,r$}
\FOR{$s = 0, \ldots, L-m-1$}
  \STATE{$W_s \la \ze{qm}{rs}f_s+\ze{qm}{r(m+s)}f_{m+s}$ }
\ENDFOR
\FOR{$s = L-m, \ldots, m-1$}
  \STATE{$W_s \la \ze{qm}{rs}f_s$ }
\ENDFOR
\STATE{$\seq{V}{\ell}{0}{m-1}\la \fft{\seq{W}{s}{0}{m-1}}$}
\RETURN{$\seq{V}{k}{0}{m-1}$}
 \end{algorithmic}
\end{algorithm}
\end{minipage}
\hfill
\begin{minipage}[t]{0.49\textwidth}
\centering
\begin{algorithm}[H]
\caption{\backwardTwo~\Bdesc{$r$}{$p=2$}{complex}}
\begin{algorithmic}\hyperlabel{alg:backwardTwo}
\INPUT{$\seq{F}{k}{0}{m-1},\param,r$}
\STATE{ $\seq{W}{s}{0}{m-1}\la \ifft{\seq{F}{\ell}{0}{m-1}}$}
\FOR{$s = 0, \ldots, m-1$}
  \STATE{$V_s \la \ze{qm}{-sr}W_s$}
\ENDFOR
\FOR{$s = m, \ldots, L-1$}
  \STATE{$V_s \la \ze{qm}{-(s-m)r}W_{s-m}$}
\ENDFOR
\RETURN{$\seq{V}{j}{0}{L-1}$}
\end{algorithmic}
\end{algorithm}
\end{minipage}
\end{multalg}

\section{Centered convolutions in one dimension}\label{sec:centeredCase}
In certain applications, it is convenient to center the data within
the input array. While it is possible to multiply the output of the
uncentered transform derived in \cref{sec:standardCase} pointwise by
the primitive roots of unity $\{\zeta^{-k\floor{qm/2}}_{qm}\}_{k=0}^{qm-1}$ to
obtain a centered transform, this
is inefficient if $qm$ is odd because it requires $qm-1$ extra
complex multiplies. When~$qm$ is even, one could in principle
account for these required signs within the multiplication
routine (with no extra reading and
writing of the Fourier transformed data); however, this would detract from
the generality of our approach.

To handle the centered case, we build the shift directly into the
transforms. Let $p,m\in\mathbb{N}$, with $p$ even, and let
$\vf=\seq{f}{j}{-pm/2}{pm/2-1}$ be a centered array
(obtained by symmetrically padding an array of length $L$ to
length $pm$, if needed). Implicit padding to length $qm$ can be
accomplished with the transform
\begin{equation}\label{eqn:centDFT}
  F_{k}= \sum_{j = -pm/2}^{pm/2-1}\ze{qm}{kj}f_{j}, \ \ \ k \in\{0, \ldots, qm - 1\}.
\end{equation}
Separating this sum and shifting the indices, we obtain
\begin{align*}
  F_{k}&=\sum_{j = 0}^{pm/2-1}\ze{qm}{kj}f_{j}+
  \sum_{j = -pm/2}^{-1}\ze{qm}{kj}f_{j}=
  \sum_{j = 0}^{pm/2-1}\ze{qm}{kj}f_{j}+
  \ze{q}{-kp/2}\sum_{j = 0}^{pm/2-1}\ze{qm}{kj}f_{j-pm/2}.
\end{align*}
Just as in \cref{sec:standardCase}, we reindex our sum (using the fact that
$p$ is even):
$$j = tm+s, \ \ t\in\{0,
\ldots,\frac{p}{2}-1\}, \ \ s\in\{0,\ldots, m-1\},$$
$$k = q\ell+r, \ \ \ell\in\{0, \ldots, m-1\}, \ \ r\in\{0,\ldots, q-1\}.$$
Then \cref{eqn:centDFT} can be computed with $q$ DFTs of size $m$:
\begin{equation}\label{eqn:forCent}
  F_{q\ell+r}=\sum_{s = 0}^{m-1}\ze{m}{\ell s}w_{r,s},
\end{equation}
where
\begin{equation}\label{eqn:wDef}
  w_{r,s}\deq\ze{qm}{r s}\left(\sum_{t=0}^{p/2-1}\ze{q}{rt}
  \left[f_{tm+s}+\ze{2q}{-rp}f_{tm+s-pm/2}\right]\right).
\end{equation}

Because the transformed data is not centered, the inverse transform is
identical to \cref{eqn:invStan}; one must only be careful to store the
output values of the inverse transform in the correct locations.
Pseudocode for a centered convolution is given in
\cref{alg:convolveCentered}. Pseudocode for centered transforms when
$p=2$ is given in \cref{alg:forwardCentered,alg:backwardCentered}. For $p > 2$, see \cref{sub:centInner}.

\begin{multalg}[htbp]
\begin{minipage}[t]{0.49\textwidth}
\centering
\begin{algorithm}[H]
\caption{\forwardCentered~\Fdesc{$r$}{$p=2$}{centered}}
\begin{algorithmic}\hyperlabel{alg:forwardCentered}
\INPUT{$\seq{f}{j}{0}{L-1},\param,r$}
\STATE{$\Hf\la \floor{L/2}$}
\FOR{$s = 0, \ldots, m - \Hf - 1$}
  \STATE{$W_{s} \la \ze{qm}{rs}{f_{\Hf+s}}$}
\ENDFOR
\FOR{$s = m - \Hf, \ldots, L -\Hf-1$}
  \STATE{$W_{s} \la \ze{qm}{r(s-m)}{f_{\Hf+s-m}} + \ze{qm}{rs}{f_{\Hf+s}}$}
\ENDFOR
\FOR{$s = L - \Hf, \ldots, m - 1$}
  \STATE{$W_{s} \la \ze{qm}{r(s-m)}{f_{\Hf+s-m}}$}
\ENDFOR
\STATE{$\seq{V}{s}{0}{m-1} \la \fft{\seq{W}{s}{0}{m-1}}$}
\RETURN{$\seq{V}{k}{0}{m-1}$}
\end{algorithmic}
\end{algorithm}
\end{minipage}
\hfill
\begin{minipage}[t]{0.49\textwidth}
\centering
\begin{algorithm}[H]
\caption{\backwardCentered~\Bdesc{$r$}{$p=2$}{centered}}
\begin{algorithmic}\hyperlabel{alg:backwardCentered}
\INPUT{$\seq{F}{k}{0}{2m-1},\param,r$}
\STATE{$\Hf\la \floor{L/2}$}
\STATE{ $\seq{W}{s}{m}{2m-1}\la \ifft{\seq{F}{s}{m}{2m-1}}$}
\FOR{$s = m-\Hf, \ldots, m-1$}
  \STATE{$V_{\Hf+s-m} \la \ze{qm}{-r(s-m)}W_{s}$}
\ENDFOR
\FOR{$s = 0, \ldots, L-\Hf-1$}
  \STATE{$V_{\Hf+s} \la \ze{qm}{-rs}W_{s}$}
\ENDFOR
\RETURN{$\seq{V}{j}{0}{L-1}$}
\end{algorithmic}
\end{algorithm}
\end{minipage}
\end{multalg}

\section{Hermitian symmetric convolutions in one dimension}
\label{sec:hermitianCase}
A centered array $\vf=\seq{f}{j}{-pm/2+1}{pm/2-1}$ is Hermitian
symmetric if $f_j = \conj{f_{-j}}$ (where the bar denotes complex conjugation)
for all $j$.
Due to the symmetry of the data, one only has to store approximately half of the
input values, as the rest of the data can be computed (by taking the conjugate)
when needed. An array is Hermitian symmetric if and only if its DFT is
real valued. Because of this, Hermitian symmetric data occurs naturally in many
applications, including pseudospectral methods for partial differential
equations.

One can use the centered transforms from \cref{sec:centeredCase} to
develop Hermitian transforms. The forward transform is given by
\cref{eqn:forCent}, where
\begin{equation}\label{eqn:wHerDef}
  w_{r,s}\deq\ze{qm}{r s}\left(\sum_{t=0}^{p/2-1}\ze{q}{rt}
  \left[f_{tm+s}+\ze{2q}{-rp}\conj{f_{pm/2-tm-s}}\right]\right).
\end{equation}
Note that we only require $f_j$ for $j\in\{0, \ldots,
pm/2\}$. Furthermore, we have the Hermitian symmetry $w_{r,s} = \conj{w_{r,-s}}$
(which holds since the DFT of $\seq{w}[r,s]{s}{0}{m-1}$ produces real-valued
output) so we only need to compute
$w_{r,s}$ for $s\in\{0,\ldots,\floor{m/2}+1\}$, using a complex-to-real
DFT to compute each residue contribution.

The inverse transform is once again given by \cref{eqn:forStan}. Here, the
key difference is that because the input is real, the output is
Hermitian symmetric, so that real-to-complex DFTs can be
used. Pseudocode for a Hermitian convolution is given in
\cref{alg:convolveCentered}. Pseudocode for the Hermitian transforms
when $p=2$ is given in \cref{alg:forwardHermitian,alg:backwardHermitian}. For $p > 2$, see \cref{sub:centInner}.

\begin{multalg}[htbp]
\begin{minipage}[t]{0.49\textwidth}
\centering
\begin{algorithm}[H]
\caption{\forwardHermitian~\Fdesc{$r$}{$p=2$}{Hermitian}}
\begin{algorithmic}\hyperlabel{alg:forwardHermitian}
\INPUT{$\seq{f}{j}{0}{\Hc},\param,r$}
\STATE{$\Hc\la \ceil{L/2}$}
\STATE{$e=\floor{m/2}+1$}
\FOR{$s = 0, \ldots, m - \Hc $}
  \STATE{$W_s \la \ze{qm}{rs}{f_s}$}
\ENDFOR
\FOR{$s = m - \Hc+1, \ldots, e-1$}
  \STATE{$W_s \la \ze{qm}{rs}{f_s} + \ze{qm}{r(s-m)}{\conj{f_{m-s}}}$}
\ENDFOR
\STATE{$\seq{V}{\ell}{0}{m-1} \la \crfft{\seq{W}{s}{0}{e-1}}$}
\RETURN{$\seq{V}{k}{0}{2m-1}$}
\end{algorithmic}
\end{algorithm}
\end{minipage}
\hfill
\begin{minipage}[t]{0.49\textwidth}
\centering
\begin{algorithm}[H]
\caption{\backwardHermitian~\Bdesc{$r$}{$p=2$}{Hermitian}}
\begin{algorithmic}\hyperlabel{alg:backwardHermitian}
\INPUT{$\seq{F}{k}{0}{2m-1},\param, r$}
\STATE{$\Hc\la \ceil{L/2}$}
\STATE{$e\la \floor{m/2}+1$}
\STATE{ $\seq{W}{s}{0}{e-1}\la \icrfft{\seq{F}{s}{0}{m-1}}$}
\FOR{$s = 0, \ldots,  m-e$}
  \STATE{$V_s\la \ze{qm}{-rs}W_s$}
\ENDFOR
\FOR{$s=m-\Hc+1, \ldots, m-e$}
  \STATE{$V_{m-s} \la \ze{qm}{-r(m-s)}\conj{W_s}$}
\ENDFOR
\IF{$m$ is even}
  \STATE{$V_{e-1} \la  \ze{2q}{-r}W_{e-1}$}
\ENDIF
\RETURN{$\seq{V}{j}{0}{\Hc}$}
\end{algorithmic}
\end{algorithm}
\end{minipage}
\end{multalg}

\section{Multidimensional convolutions}\label{sec:multidim}

\begin{figure}[htbp]
\begin{center}
\includegraphics[width=\linewidth]{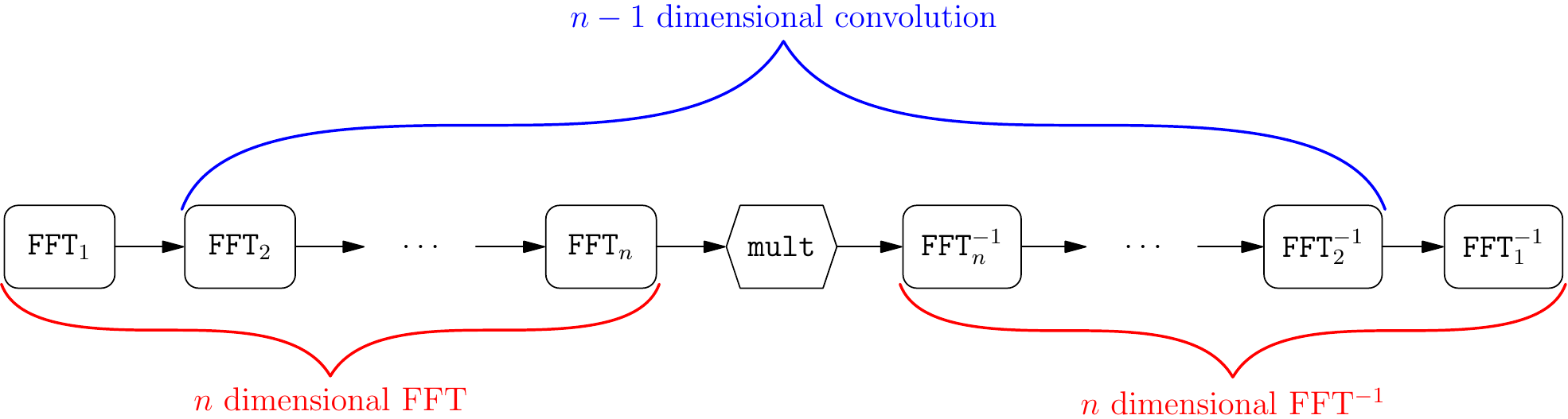}
\caption{Recursive computation of an $n$-dimensional convolution.}
\label{fig:multiDimConv}
\end{center}
\end{figure}

An $n$-dimensional convolution is conventionally computed by performing an FFT
of size $N_1\times\ldots\times{}N_n$, applying the specified multiplication
operator on the transformed data, and then performing an inverse FFT back to the
original space. However, as described in~\cite{Bowman11}
and~\cite{Roberts18}, a better alternative is to decompose the
$n$-dimensional convolution recursively into $\prod_{i=2}^{n}N_i$
FFTs in the first dimension, followed by $N_1$ convolutions of
dimension $n-1$, and finally $\prod_{i=2}^{n}N_i$ inverse FFTs in the first dimension.
This is illustrated in~\cref{fig:multiDimConv}. At the innermost level,
a recursive multidimensional convolution
reduces to a one-dimensional convolution.

The most important advantage of decomposing a multidimensional
convolution is that one can reuse the work buffer for each
subconvolution, reducing the total memory footprint. These storage
savings are attainable regardless of whether explicit or implicit
dealiasing is used for the underlying padded FFTs.

\begin{figure}[htbp]
  \begin{center}
    \includegraphics{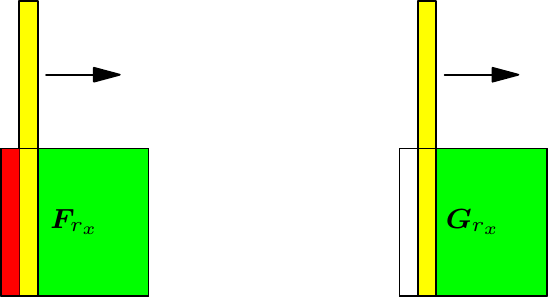}
    \caption{The reuse of memory to compute the contribution of a
      single $x$ residue to a 2D binary convolution with two inputs
      and one output: a 1D padded~$y$ FFT is applied to columns of
      $\vF_{r_x}$ and~$\vG_{r_x}$ to produce the two stacked yellow
      columns that are fed
      to the multiplication  operator, producing one stacked column to
      be inverse~$y$ transformed into a single column (like
      the red one shown on the left). The upper column
      is then reused for processing subsequent columns.
      }\label{conv2psimp4}
    \end{center}
\end{figure}

For example, the memory management for a single-threaded 2D
padded complex convolution for $A=2$ and $B=1$ is shown in
\cref{conv2psimp4}. For each \hbox{$r_x\in \{0,1,\ldots, q_x-1\}$}, the residue
contribution to the padded $x$ FFT of the input buffers is stored in the square boxes.
A padded FFT of each input is then performed
in the $y$ direction, column-by-column, using a one-dimensional work buffer,
to produce a single column of the Fourier-transformed image,
depicted in yellow.
The Fourier transformed columns of two inputs~$F$ and~$G$ are then multiplied
pointwise and stored back into the $F$~column. At this point, the
inverse $y$ transform can then be performed, with the truncated result
stored in the
lower half of the column, next to the previously processed data shown in
red. This process is repeated on the remaining columns, shifting and
reusing the work buffers. Once all the columns have been processed, an
inverse transform in the~$x$ direction produces the final~$r_x$ contribution
to the convolution.

The reuse of subconvolution work memory allows the convolution to be
computed using less total memory: for a $d$-dimensional $p/q$ padded
convolution, the work memory requirement is $(A+B)pm L^{d-1}$ complex
words, where $p=\ceil{L/m}$ (not counting the storage requirements for
the input data).
In contrast, explicit padding requires a typically much larger buffer
of size $E\doteq C(q/p)^dL^d$ words, where $C=\max(A,B)$.
For example, computing a $d$-dimensional dealiased convolution
implicitly for $A=2$ and $B=1$ with padding ratio $p/q=1/2$ and
$m=L \gg 1$ asymptotically requires a storage of $3E/2^{d+1}$.
In particular, in one dimension, the general formulation of implicit
padding requires $3/4$ of the work memory required by explicit padding.
For a $2/3$ padding ratio, implicit padding requires a storage of
$\left(2/3\right)^{d-1} E$.
In addition to having reduced memory requirements,
a dealiased multidimensional convolution decomposed in this way is
significantly faster than a conventional implementation due to better
data locality and the elimination of transforms of
data known \emph{a priori} to be zero.

For Hermitian-symmetric data, we assume that the origin is in the
center of the unsymmetrized domain (only about half of which is
retained). The outer convolutions are centered, while the innermost
convolution is Hermitian. The input is assumed to be Hermitian
symmetric on the hyperplane orthogonal to the innermost direction.

The FFTs in multithreaded convolutions
can be parallelized by dividing the $\prod_{i=2}^nN_i$
one-dimensional FFTs among several threads. Similarly, the $N_1$ subconvolutions
can be parallelized over $T\le N_1$ threads using $T$ work buffers
\cite{Roberts18}.

\section{Numerical implementation}\label{sec:optimization}

In this section, we describe several optimizations that significantly
improve the performance of the underlying one-dimensional padded/unpadded FFTs.

\subsection{Inner loop optimization}\label{sub:inner}
Consider the complex case described in \cref{sec:standardCase}.
The preprocessing done in \cref{eqn:forStan} is itself a padded
FFT from size $p$ to size~$q$ and the post-processing
in~\cref{eqn:invStan} is an unpadded FFT from size $q$ to size $p$.
Thus, if $p$ and $q$ share a common factor, one can use these equations
recursively.

We redefine $q$ in \cref{eqn:pqdef1} as
the smallest positive multiple of $p$ such that $qm \geq M$:
\begin{equation*}
  n\deq \ceil{\frac{M}{pm}}, \ \ q\deq np.
\end{equation*}
To compute the preprocessing in \cref{eqn:forStan}, we let $r = un +
v$, where $u \in \{0, \ldots, p-1\}$, and $v \in\{0, \ldots, n-1\}$. Then the forward
transform becomes
\begin{equation*}
  F_{q\ell+un+v} =\sum_{s=0}^{m-1}\ze{m}{\ell s}\ze{qm}{(un+v)s}
  \sum_{t=0}^{p-1}\ze{p}{ut}\left(\ze{q}{vt}f_{tm+s}\right),
\end{equation*}
so the sum over $t$ can be computed using $n$ DFTs of size $p$. Similarly, the
post-processing in \cref{eqn:invStan},
\begin{equation}\label{eqn:invInner}
  f_{tm+s} = \frac{1}{qm}\sum_{v=0}^{n-1}\ze{np}{-tv}
  \sum_{u=0}^{p-1}\ze{p}{-tu}\ze{qm}{-s(un+v)}
  \sum_{\ell=0}^{m-1}\ze{m}{-s\ell}F_{q\ell+un+v},
\end{equation}
requires the sum of $n$ DFTs of size~$p$. In this optimization we consider each
$v$ to be a residue; the full transform then has $n$ residue contributions of
size $pm$. Note that the sum over $v$ is \emph{not} a DFT as the input depends on both $v$ and~$t$.
Pseudocode for the inner loop is given in \cref{alg:forwardInner,alg:backwardInner}.
\begin{multalg}[htbp]
\begin{minipage}[t]{0.49\textwidth}
\centering
\begin{algorithm}[H]
\caption{\forwardInner~\Fdesc{$v$}{$p>2$}{complex}}
\begin{algorithmic}\hyperlabel{alg:forwardInner}
\INPUT{$\seq{f}{j}{0}{L-1},\param,v$}
\STATE{$p\la \ceil{L/m}$}
\FOR{$t = 0, \ldots, p-2$}
  \FOR{$s = 0, \ldots, m-1$}
    \STATE{$W_{tm+s} \la \ze{q}{vt}f_{tm+s}$}
  \ENDFOR
\ENDFOR
\STATE{$t_0\la p-1$}
\FOR{$s = 0, \ldots, L-t_0m-1$}
  \STATE{$W_{t_0m+s}\la  \ze{q}{vt_0}f_{t_0m+s}$}
\ENDFOR
\FOR{$s = L-t_0m, \ldots, m-1$}
  \STATE{$W_{t_0m+s}\la 0$}
\ENDFOR
\FOR{$s = 0, \ldots, m-1$}
  \STATE{$\seq{W}[um+s]{u}{0}{p-1} \la \fft{\seq{W}[tm+s]{t}{0}{p-1}}$}
\ENDFOR
\FOR{$u=0,\ldots,p-1$}
  \FOR{$s=0,\ldots,m-1$}
    \STATE{$W_{um+s} \la \ze{qm}{(un+v)s}W_{um+s}$}
  \ENDFOR
  \STATE{$\seq{F}[um+\ell]{\ell}{0}{m-1} \la \fft{\seq{W}[um+s]{s}{0}{m-1}}$}
\ENDFOR
\RETURN{$\seq{F}{k}{0}{pm-1}$}
\end{algorithmic}
\end{algorithm}
\end{minipage}
\hfill
\begin{minipage}[t]{0.49\textwidth}
\centering
\begin{algorithm}[H]
\caption{\backwardInner~\Bdesc{$v$}{$p>2$}{complex}}
\begin{algorithmic}\hyperlabel{alg:backwardInner}
\INPUT{$\seq{F}{k}{0}{pm-1},\param,v$}
\STATE{$p\la \ceil{L/m}$}
\FOR{$u=0,\ldots,p-1$}
  \STATE{$\seq{W}[um + s]{s}{0}{m-1} \la \ifft{\seq{F}[um+\ell]{\ell}{0}{m-1}}$}
\ENDFOR
\FOR{$u = 0, \ldots, p-1$}
  \FOR{$s = 0, \ldots, m-1$}
    \STATE{$W_{um+s} \la \ze{qm}{-s(un+v)}W_{um+s}$}
  \ENDFOR
\ENDFOR
\FOR{$s = 0, \ldots, m-1$}
  \STATE{$\seq{W}[tm+s]{t}{0}{p-1} \la \ifft{\seq{W}[um+s]{u}{0}{p-1}}$}
\ENDFOR
\FOR{$t=0,\ldots, p-2$}
  \FOR{$s=0,\ldots, m-1$}
    \STATE{$V_{tm+s}\la \ze{q}{-tv}W_{tm+s}$}
  \ENDFOR
\ENDFOR
\STATE{$t_0\la p-1$}
\FOR{$s = 0, \ldots, L-mt_0-1$}
  \STATE{$V_{t_0m+s}\la \ze{q}{-t_0v}W_{t_0m+s}$}
\ENDFOR
\RETURN{$\seq{V}{j}{0}{L-1}$}
\end{algorithmic}
\end{algorithm}
\end{minipage}
\end{multalg}

\subsection{Inner loop for centered and Hermitian arrays}\label{sub:centInner}
Just as we did for the complex (uncentered) case, we can apply the same recursive
techniques to the centered case described in \cref{sec:centeredCase}. The summation in \cref{eqn:wDef} is itself a padded DFT from size~$p/2$ to size $q$.
Therefore, if $q$ shares a factor with $p/2$, we can use the transforms in
\cref{sub:inner} to preform the preprocessing.

In our implementation, we consider two cases. If $p=2$, we
directly sum the two terms in $w_{r,s}$. If $p > 2$
(with $p$ assumed to be even), we define
\begin{equation*}
  n\deq \ceil{\frac{2M}{pm}},\ \ q\deq n\frac{p}{2}.
\end{equation*}
Then, letting $r = un+v$, where $u \in \{0, \ldots, p/2-1\}$, and
$v \in \{0, \ldots, n-1\}$, we compute
\begin{equation*}
  w_{un+v,s}=\ze{qm}{(un+v) s}\sum_{t=0}^{p/2-1}\ze{p/2}{ut}
  \left[\ze{q}{vt}\left(f_{tm+s}+\ze{n}{-v}f_{tm+s-pm/2}\right)\right].
\end{equation*}
For each value of $v$, each of these sums is a DFT of length $p/2$. Then, using
\cref{eqn:forCent}, we compute
\begin{equation}\label{eqn:forCentIn}
  F_{q\ell+un+v}=\sum_{s = 0}^{m-1}\ze{m}{\ell s}w_{un+v,s}.
\end{equation}
The inverse transform is the same as \cref{eqn:invInner}. Pseudocode for the
centered case is given in \cref{alg:forwardCenteredInner,alg:backwardCenteredInner}.

These equations also apply to the Hermitian case
(\cref{sec:hermitianCase}), using $f_{j}=\conj{f_{-j}}$ whenever $j<0$; however,
unlike the outer FFTs of length $m$,
the preprocessing and postprocessing stages use complex-to-complex FFTs.
Pseudocode for the Hermitian case is given in
\cref{alg:forwardHermitianInner,alg:backwardHermitianInner}.

\subsection{Conjugate symmetry optimization}\label{sub:conjSymOpt}
We can also exploit conjugate symmetries in the
primitive roots. This has been used in previous
implementations of implicit dealiasing \cite{Bowman11} for centered
convolutions, as discussed in \cref{sec:centeredCase}. Here we extend
the technique to more general situations, including $p > 2$.

First, we consider the complex case (\cref{sec:standardCase}). We use the
inner-loop formulation from \cref{sub:inner} for generality.
Assuming that $v\notin\{0, n/2\}$, define
\begin{equation*}
A_{s,t,v}\deq\ze{q}{vt}\text{Re }f_{tm+s}, \quad
B_{s,t,v}\deq i\ze{q}{vt}\text{Im }f_{tm+s}.
\end{equation*}
Then we have
\begin{equation*}
F_{q\ell+un+v}=\sum_{s=0}^{m-1}\ze{m}{\ell s}\ze{qm}{(un+v)s}
\left[\sum_{t = 0}^{p-1}\ze{p}{ut}\left(A_{s,t,v}+B_{s,t,v}\right)\right],
\end{equation*}
\begin{equation*}
F_{q\ell+un-v}=\sum_{s=0}^{m-1}\ze{m}{\ell s}\ze{qm}{(un-v)s}
\left[\sum_{t = 0}^{p-1}\ze{p}{ut}{\left(\conj{A_{s,t,v}}-
\conj{B_{s,t,v}}\right)}\right].
\end{equation*}
This allows us to compute $F_{q\ell+un+ v}$ and $F_{q\ell+un-v}$
together efficiently.

A similar optimization in the centered and Hermitian
cases (described in \cref{sec:centeredCase,sec:hermitianCase}) is obtained with
\begin{equation*}
A_{s,t,v}\deq\ze{q}{vt}\left(\text{Re }f_{tm+s}+
\ze{n}{-v}\text{Re }f_{tm+s-pm/2}\right),
\end{equation*}
\begin{equation*}
B_{s,t,v}\deq i\ze{q}{vt}\left(\text{Im }f_{tm+s}+
\ze{n}{-v}\text{Im }f_{tm+s-pm/2}\right),
\end{equation*}
which allows us to compute \cref{eqn:forCentIn} using
\begin{equation*}
w_{un+v,s}=\ze{qm}{(un+v)s}\left[\sum_{t = 0}^{p/2-1}\ze{p/2}{ut}
\left(A_{s,t,v}+B_{s,t,v}\right)\right],
\end{equation*}
\begin{equation*}
w_{un-v,s}=\ze{qm}{(un-v)s}\left[\sum_{t = 0}^{p/2-1}\ze{p/2}{ut}
{\left(\conj{A_{s,t,v}}-\conj{B_{s,t,v}}\right)}\right].
\end{equation*}

\subsection{Overwrite optimization}
For certain padded FFTs, where all residues are computed at once,
the input array is large enough to hold all but one
of the residue contributions. We have designed specialized algorithms
for such cases, with one residue contribution written to the output buffer and
the others stored in the input buffer.

The overwrite optimization is particularly advantageous in the complex
case (\cref{sec:standardCase}) when
$M \le 2L$ and $m=L$, so that $p=1$ and $q=2$.
In this case, the input buffer already contains the preprocessed data
for $r=0$ and the preprocessed data $\ze{qm}{s} f_s$ for $r=1$ is
written to the output buffer.
The input and output buffers are then individually Fourier transformed
to obtain the required residues. The backwards transform does the
reverse operation: it first performs inverse Fourier transforms on the
input and output buffers, then adds products of the output buffer
and roots of unity to the input buffer.

To use the overwrite optimization for computing a one-dimensional convolution,
the number of inputs $A$ must be at least as large as the number of
outputs~$B$. Under this same restriction, the overwrite optimization
can be implemented for each dimension of a multidimensional convolution.

\subsection{Loop optimizations}
If the overwrite optimization is not applicable, other
data flow improvements may be possible.
Suppose that we compute a block of~$D$ residues at a time as suggested in
\cref{sub:oneRes}. Normally, the contribution to the inverse
padded Fourier transform from the block containing residue $0$ is stored in an
accumulation buffer; the contributions from the remaining
residues are then added to this buffer by iterating over the other residue
blocks. If there are no other residue blocks, a separate accumulation
buffer is not needed; one can accumulate the residues
entirely within the input buffer.

Another optimization is possible when $A > B$ and there are exactly two
residue blocks, which we label $0$ and $r$. In this case, we compute
all $A$ forward-padded FFT contributions to residue block $0$
in an output buffer $F$ and apply the multiplication operator, freeing
up the storage in $F$ associated with $A-B$ inputs. We then
transform the contributions to residue block $r$, $A-B$ inputs at a
time, each time writing $A-B$ inverse-transformed contributions
from residue block $0$ to the input buffer. Once all $A$
contributions to residue block $r$ have been forward transformed, we
apply the multiplication operator and accumulate the contributions
from the inverse transform in the input buffer.

\section{Numerical results}\label{sec:results}
This work presents several algorithms for computing padded FFTs.
Which algorithm is optimal for a given problem depends on
the number of inputs $A$ and outputs $B$, the multiplication operator \mult,
the input data length~$L$, the padded length $M$, the number of copies $C$
of the transform to be computed simultaneously, and the stride $S$
between successive data elements of each copy. We only
consider the efficiently packed case, where the distance in memory between the
first elements of each copy is one.
Vectorized and parallelized C++ versions of
these algorithms have been implemented in the open-source library
{\tt FFTW++}~\cite{fftwpp}.

We determine the fastest algorithm for a given problem empirically,
scanning over the underlying FFT size $m$, the number $D$ of residues
to be computed at a time, and whether or not to use in-place or
out-of-place FFTs. These parameters then determine the values of $p$
and~$q$ to use in our padded FFT algorithms.

Our optimizer measures the time required to compute a one-dimensional
in-place dealiased convolution for a particular set of parameters, using
the given multiplication routine.
As described in \cref{sec:multidim}, multidimensional convolutions
are decomposed into a sequence of padded FFTs, a one-dimensional
convolution, and then a sequence of inverse padded FFTs. We assume that
the padded/unpadded FFT pairs can be optimized independently in
each dimension. This is accomplished by performing one-dimensional
convolutions using each padded FFT pair, without calling the
multiplication routine. In practice, this decoupling works well and
leads to efficient optimization of multidimensional convolutions.
If the outermost pair FFT${}_1$ and FFT${}_1^{-1}$
in~\cref{fig:multiDimConv} are to be multithreaded over $T$ threads,
optimization of the remaining FFTs should be performed over $T$
concurrent copies, to simulate the execution environment.

We benchmarked our algorithms with a liquid-cooled Intel i9-12900K
processor (5.2GHz, 8 performance cores) on an ASUS ROG Strix Z690-F
motherboard with 128GB of DDR5 memory (5GHz), using
version  12.2.1 of the {\tt GCC} compiler with the optimizations
{\tt -Ofast -fomit-frame-pointer -fstrict-aliasing -ffast-math}.
The underlying FFTs were computed with version 3.3.10 of the adaptive
{\tt FFTW} \cite{fftw,Frigo05} library under the Fedora 37 operating
system. Multithreading was implemented with the {\tt OpenMP} library.

\subsection{One-dimensional convolutions}
In~\cref{fig:timings1-T1} we plot median execution times (normalized to
$L\log_2 L$) for one-dimensional in-place
convolutions of $L$ complex words (with $M=2L$) over one thread
for explicit zero padding using in-place (IP) or out-of-place
(OP) FFTs, implicit dealiasing \cite{Roberts18}, and hybrid dealiasing.
We see that hybrid dealiasing is much faster than both explicit and
implicit dealiasing for large sizes and is generally competitive with
optimized out-of-place explicit algorithms (which it reduces to) for
small sizes. In~\cref{fig:timings1-T8} we plot the normalized times for
the same one-dimensional convolutions parallelized over 8 threads.
In these plots, only power-of-two sizes are benchmarked since these are
optimal FFT sizes.

\begin{figure}[tbhp]
\begin{minipage}{0.49\linewidth}
\Figure{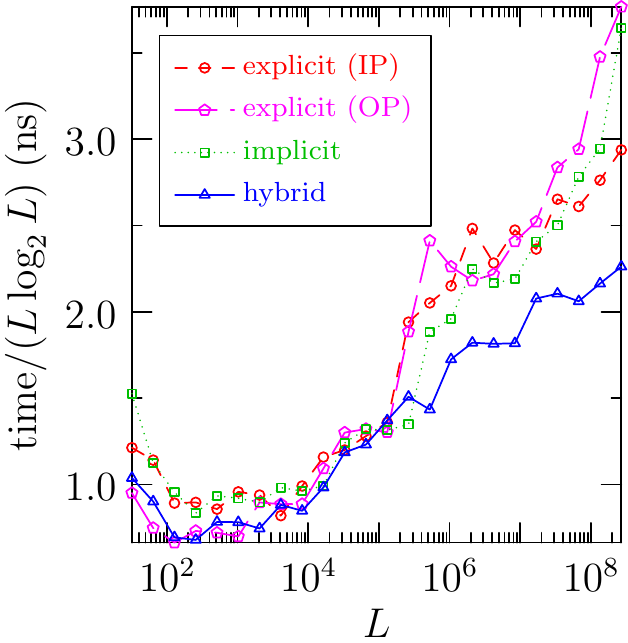}{In-place 1D complex convolutions of length $L$
  with $A=2$ and $B=1$ on 1 thread.}
\end{minipage}
\,
\begin{minipage}{0.49\linewidth}
\Figure{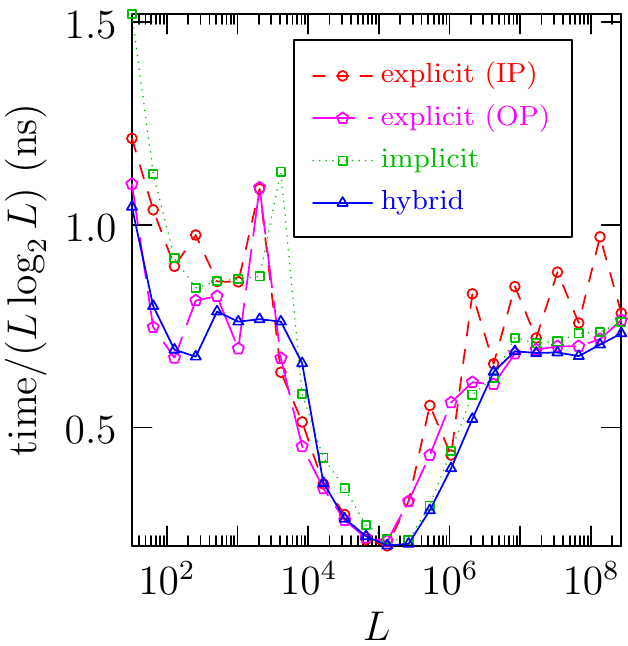}{In-place 1D complex convolutions of length $L$
  with $A=2$ and $B=1$ on 8 threads.}
\end{minipage}
\end{figure}

In~\cref{fig:timingsh1-T1,fig:timingsh1-T8} we plot the
normalized times for one-dimensional Hermitian convolutions of size $L$
padded to $M=3L/2$ over 1 thread and 8 threads, respectively.
We benchmark hybrid dealiasing separately for those sizes that are optimal for
explicit dealiasing and for implicit dealiasing.
For the implicit dealiasing algorithms
developed in \cite{Bowman11,Roberts18}, the optimal sizes are one less
than a power of two. In our implementation of hybrid dealiasing,
we normally adjust these sizes to exact powers of two to allow us to use
the overwrite optimization. For explicit dealiasing, the optimal
values of $L$ are $2\floor{\frac{2^n+2}{3}}-1$ for positive integers $n$.
We observe that hybrid dealiasing outperforms implicit dealiasing at
optimal implicit sizes and performs about as well
as explicit dealiasing at optimal explicit sizes.

\begin{figure}[tbhp]
\begin{minipage}{0.49\linewidth}
\Figure{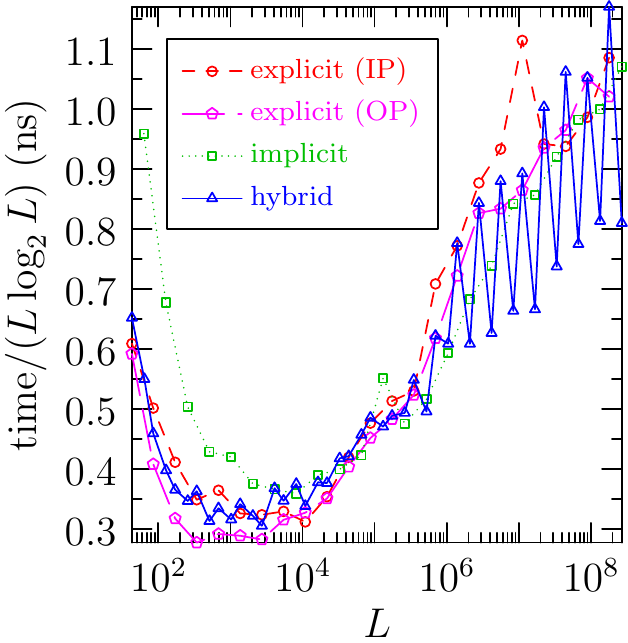}{In-place 1D Hermitian convolutions of length $L$
  with $A=2$ and $B=1$ on 1 thread.}
\end{minipage}
\,
\begin{minipage}{0.49\linewidth}
\Figure{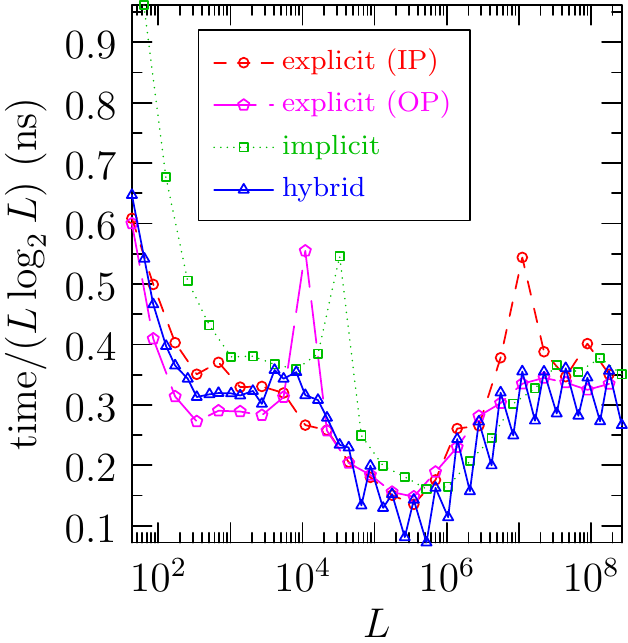}{In-place 1D Hermitian convolutions of length $L$
  with $A=2$ and $B=1$ on 8 threads.}
\end{minipage}
\end{figure}

\subsection{Two-dimensional convolutions}

In~\cref{fig:timings2-T1,fig:timings2-T8} we plot the
normalized median times for 2D complex convolutions of size
$L\times L$. On a single thread, these benchmarks show that at all sizes,
hybrid dealiasing is faster than implicit dealiasing and much faster than
explicit dealiasing. On 8 threads, hybrid dealiasing outperforms both methods except at
$L=64$. In both cases an $x$ stride of $L+2$ was used.

\begin{figure}[tbhp]
\begin{minipage}{0.49\linewidth}
\Figure{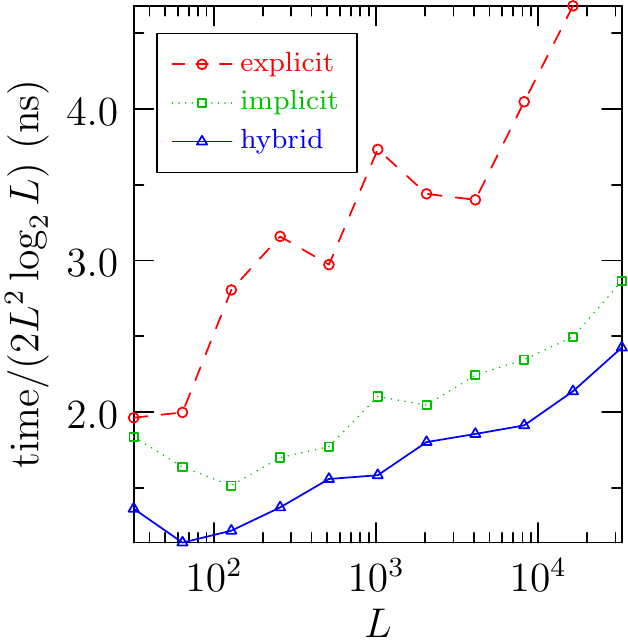}{In-place 2D complex convolutions of size $L\times L$
  with $A=2$ and $B=1$ on 1 thread.}
\end{minipage}
\,
\begin{minipage}{0.49\linewidth}
\Figure{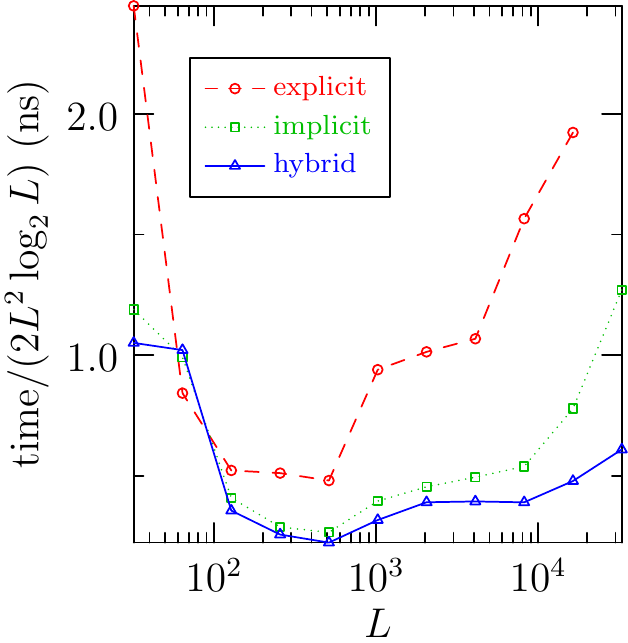}{In-place 2D complex convolutions of size $L\times L$
  with $A=2$ and $B=1$ on 8 threads.}
\end{minipage}
\end{figure}

In \cref{fig:timingsh2-T1,fig:timingsh2-T8},
which record normalized median times for two-dimensional Hermitian
convolutions of size $L\times L$, we see that hybrid dealiasing
outperforms both implicit and explicit dealiasing for all sizes.
In both cases we used an $x$ stride of~$\ceil{\frac{L}{2}}+2$.

To illustrate the sometimes counter-intuitive optimal parameters,
in the single-threaded case with $L=2048$, the optimizer chose
the parameters $m=16$, $p=128$, $q=192$, and $D=1$ in the $x$ direction and
performed the transform in place (using the inner loop optimization), whereas
in the $y$ direction, out-of-place explicit dealiasing ($m=3072$ and
$p=q=D=1$) was optimal.

\begin{figure}[tbhp]
\begin{minipage}{0.49\linewidth}
\Figure{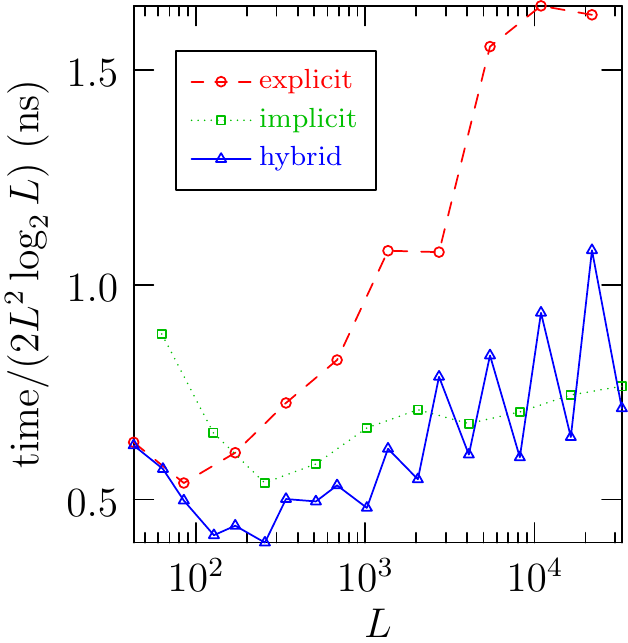}{In-place 2D Hermitian convolutions of size $L\times L$
  with $A=2$ and $B=1$ on 1 thread.}
\end{minipage}
\,
\begin{minipage}{0.49\linewidth}
\Figure{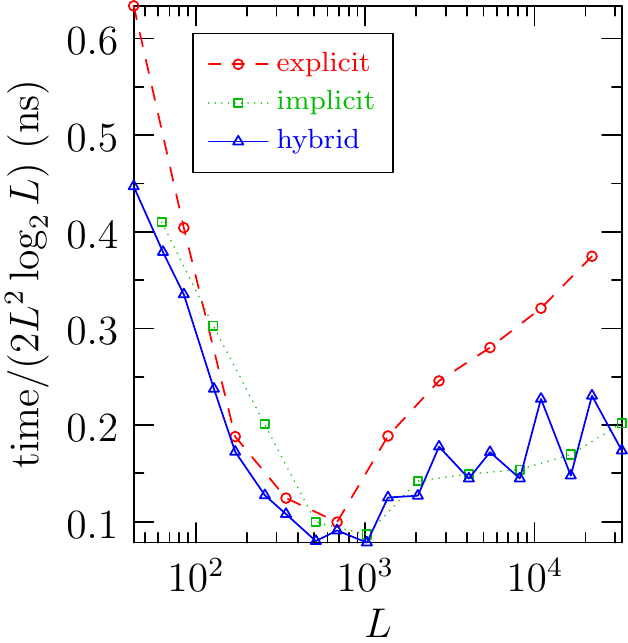}{In-place 2D Hermitian convolutions of size $L\times L$
  with $A=2$ and $B=1$ on 8 threads.}
\end{minipage}
\end{figure}

\subsection{Three-dimensional convolutions}

In~\cref{fig:timings3-T1,fig:timings3-T8}
we benchmark three-dimensional complex convolutions.
In the single-threaded case we used the $y$ stride $S_y=L+2$ and the
$x$ stride $S_yL+2$; in the multithreaded-threaded case, we used the
$y$ stride $S_y=L$ and the $x$ stride $S_yL+4$. Again, we observe that
hybrid dealiasing outperforms the other two methods.

\begin{figure}[tbhp]
\begin{minipage}{0.49\linewidth}
  \Figure{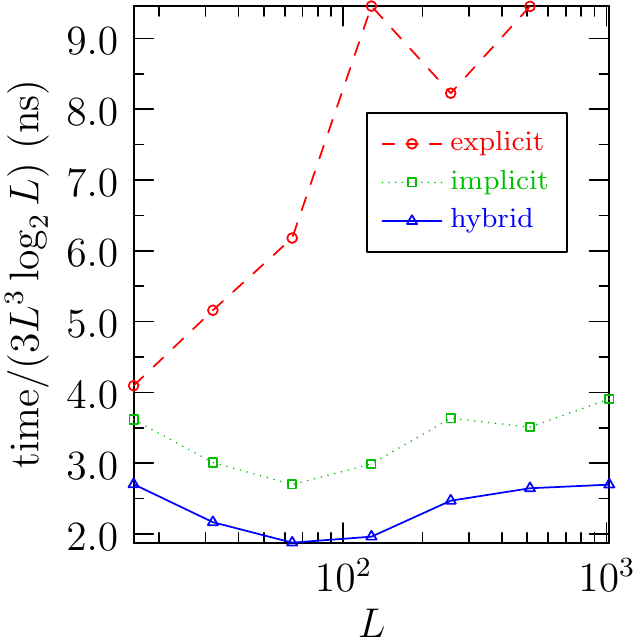}{In-place 3D complex convolutions of
    size $L\times L \times L$ with $A=2$ and $B=1$ on 1 thread.}
\end{minipage}
\,
\begin{minipage}{0.49\linewidth}
\Figure{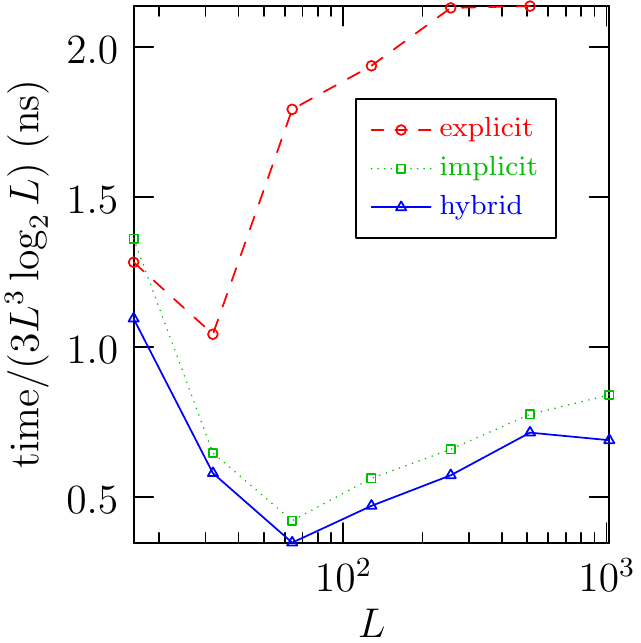}{In-place 3D complex convolutions of
    size $L\times L \times L$ with $A=2$ and $B=1$ on 8 threads.}
\end{minipage}
\end{figure}

In~\cref{fig:timingsh3-T1} we observe that for three-dimensional Hermitian
convolutions on a single thread, hybrid dealiasing performs better than
implicit dealiasing and much better than explicit dealiasing, except at
$L=21$. As seen in~\cref{fig:timingsh3-T8}, when run on 8 threads, hybrid
dealiasing is much faster than the other methods for all sizes.
In the single-threaded case we used the $y$ stride $S_y=
\ceil{\frac{L}{2}}+2$ and the $x$ stride $S_yL+2$; in the
multithreaded-threaded case, we used the $y$ stride $S_y=L$ and
the $x$ stride~$S_yL+4$.

\begin{figure}[tbhp]
\begin{minipage}{0.49\linewidth}
\Figure{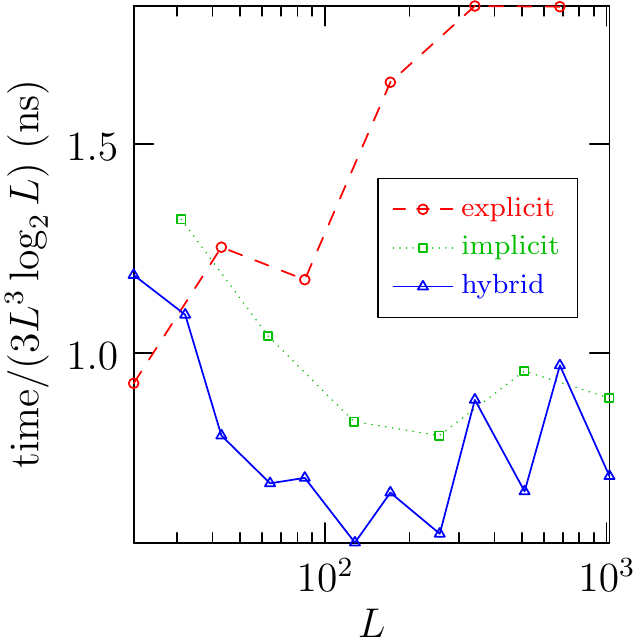}{In-place 3D Hermitian convolutions of
    size $L\times L \times L$ with $A=2$ and $B=1$ on 1 thread.}
\end{minipage}
\,
\begin{minipage}{0.49\linewidth}
\Figure{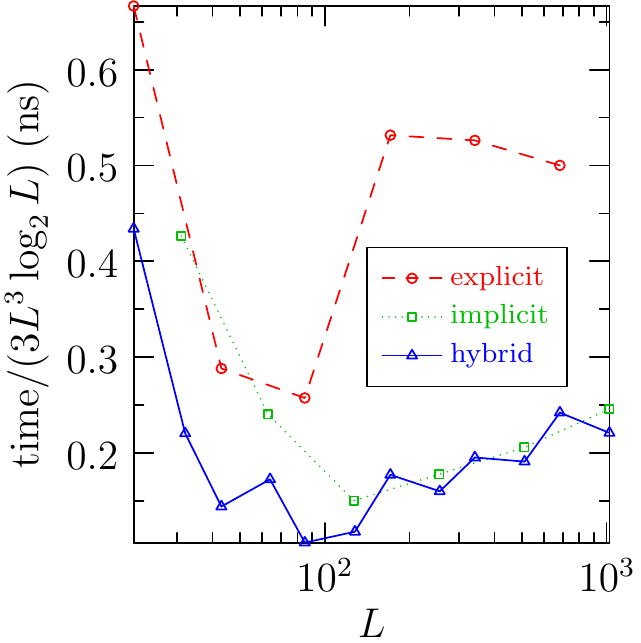}{In-place 3D Hermitian convolutions of
    size $L\times L \times L$ with $A=2$ and $B=1$ on 8 threads.}
\end{minipage}
\end{figure}

In~\cref{fig:timings3-T1I1}, we emphasize that, unlike explicit and implicit
dealiasing, hybrid dealiasing performs consistently well over a range of
arbitrary sizes.

\begin{figure}[tbhp]
\Figure{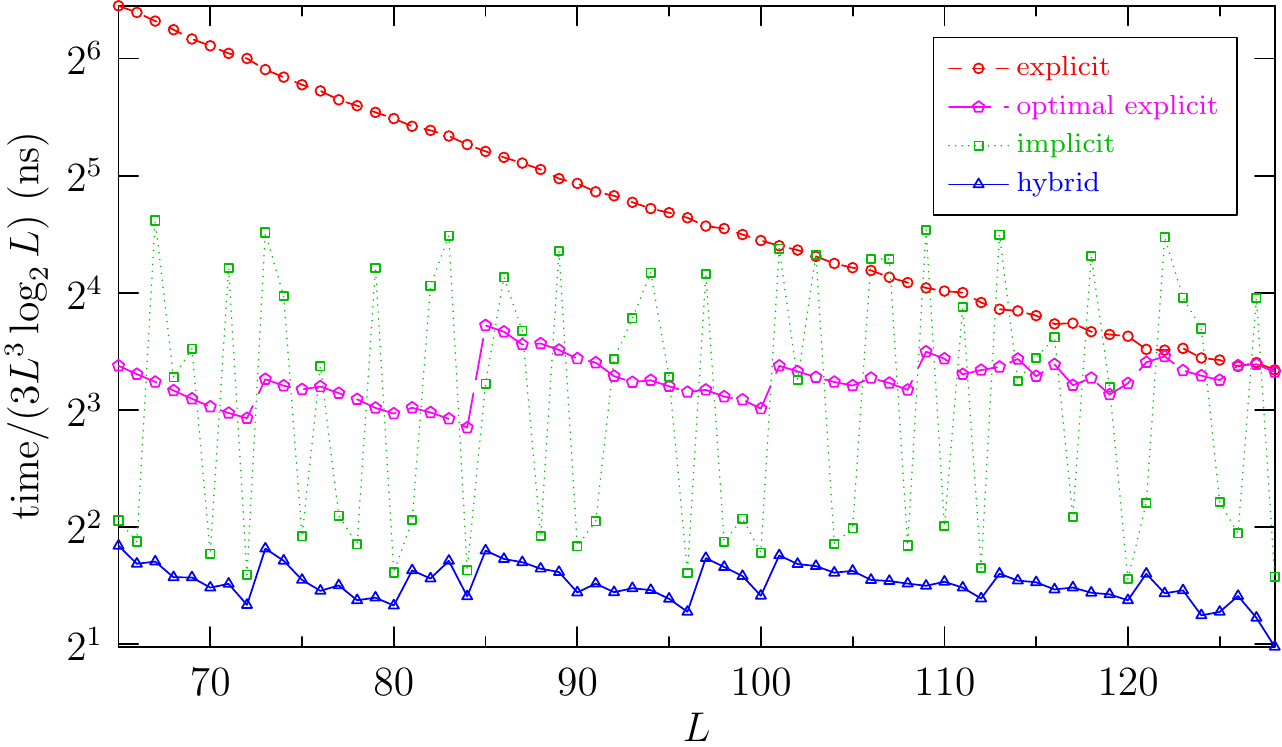}{Semi-log plot of 3D in-place complex
  convolution timings for incremental
  sizes $L\times L\times L$ with $A=2$ and $B=1$ on 1 thread. As is common
  practice, explicit dealiasing pads past $2L$ up
  to the next power of two (in this case 256). Optimal explicit dealiasing
  zero pads beyond $2L$ to the empirically determined optimal size.}
\end{figure}

\section{Future work}\label{sec:applications}
Many applications rely on real convolutions: signal
processing, image processing, and convolutional neural networks.
FFT-based convolutions are widely used in signal and image
processing and are beginning to be used to accelerate the training of
convolutional neural networks \cite{MathieuEtAl14,HighlanderRodriguez16,
LinYao19,ChitsazEtAl20}.
The algorithms described in this work operate on complex data and are
not optimal for the real input data used in such applications.
Moreover, many of these applications involve convolving input data
with a relatively small kernel. Because kernel sizes are typically
much smaller than the arrays they are convolved with, the majority of
the padded kernel is made up of zeros. In future work, we will
specialize hybrid dealiasing to real arrays of different lengths
and develop a general toolkit for diverse applications.

One of the strengths of our formulation for computing padded/unpadded
DFTs is that we decompose them into full DFTs, facilitating
the use of existing FFT algorithms.
Hybrid dealiasing does not attempt to compete with, but rather harness, modern FFT libraries.

Unfortunately, the real-to-complex case is not as simple as the complex and
Hermitian cases described in this work. For Hermitian symmetric data
(\cref{sec:hermitianCase})
each residue contribution is real, as it is part of a larger real array and
can be computed using a complex-to-real DFT.
For real data, the contribution of any residue other than zero is not, in
general, Hermitian symmetric. Therefore, one cannot simply compute the residue
contribution using a real-to-complex DFT.
While one could pack the real data into a smaller complex array and use the
algorithms described in this paper, such schemes are less efficient than
algorithms designed with real data in mind \cite{SorensenEtAl87}.

\section{Conclusion}\label{sec:conclusion}

This work combines several ideas to implement efficient dealiased
convolutions for arbitrary padding ratios. Hybrid dealiasing combines
the techniques of implicit and explicit dealiasing to exploit
optimal transform sizes of the underlying FFT library.
By exploiting smaller FFT sizes than used in \cite{Bowman11},
hybrid dealiasing can store
the required roots of unity in short unfactorized tables.
The recognition that the preprocessing in \cref{eqn:forStan} is itself
a padded FFT is crucial to an efficient implementation for general
padding ratios.

Hybrid dealiasing is ideally suited for implementing pseudospectral
simulations of partial differential equations.
For example, the three-dimensional incompressible Navier--Stokes equations
can be implemented with the Basdevant formulation using $A=3$ and $B=5$,
while the two-dimensional formulation uses $A=B=2$ \cite{Basdevant83}.
Likewise, the incompressible three-dimensional magnetohydrodynamic
equations can be implemented with $A=6$ and $B=8$.\footnote{Here we
correct an error in \cite{Roberts18}; only $A+B=14$ FFT calls are required.}
The generality of the formulation described in this paper can also be
applied to the cascade direction of $n$th-order Casimir invariants
of two-dimensional turbulence, such as
$\sum_j \omega^n(x_j)$, where $\omega$ is the scalar vorticity
 \cite{Bowman13casimir}; this requires a padding ratio of $2/(n+1)$.

The naive way of exploiting the convolution theorem is to take the
Fourier transform
of the entire input data, perform the necessary multiplication, and
then take the inverse transform.
However, if a convolution is all that is needed,
we have seen that better performance can be achieved by localizing the
computations; one can
compute the contribution to the convolution one residue at a time.
Furthermore, multidimensional convolutions can be done more efficiently
by decomposing them into an outer FFT, a lower-dimensional
convolution, and an inverse FFT. The possibility of reusing work
memory in this recursive formulation is responsible for most of the
dramatic performance gains that are observed in two and three dimensions.

\section*{Acknowledgment}
Financial support for this work was provided by grants
RES0043585 and RES0046040 from the Natural Sciences and Engineering Research
Council of Canada.

\appendix
\section{Centered and Hermitian inner loop pseudocode}

Pseudocode for computing centered padded/unpadded FFTs with the
inner loop optimization (\cref{sub:centInner})
is presented in \cref{alg:forwardCenteredInner,alg:backwardCenteredInner}.
The Hermitian case is documented in
\cref{alg:forwardHermitianInner,alg:backwardHermitianInner}.

\begin{multalg}[htbp]
\begin{minipage}[t]{0.49\textwidth}
\centering
\begin{algorithm}[H]
\caption{\forwardCenteredInner~\Fdesc{$v$}{$p>2$}{centered complex}~\deltadesc}
\begin{algorithmic}\hyperlabel{alg:forwardCenteredInner}
\INPUT{$\seq{f}{j}{0}{L-1},\param,v$}
\STATE{$\Hf\la \floor{L/2}$}
\STATE{$\ph \la \ceil{L/(2m)}$}
\STATE{$m_0\la \phm-\Hf$}
\STATE{$m_1\la L-\Hf-(\ph-1)m$}
\FOR{$s = 0, \ldots, m_0-1$}
  \STATE{$W_{s} \la f_{\Hf+s}$}
\ENDFOR
\FOR{$t = 0, \ldots, \ph-1$}
  \STATE{$m_2 \la (m_1-m)\delta_{t-(\ph-1)}+m$}
  \FOR{$s = m_0\delta_t, \ldots, m_2-1$}
    \STATE{$W_{tm+s} \la \ze{q}{v(t-\ph)}f_{(t-\ph)m+\Hf+s}+$}
    \STATE{\phantom{$W_{tm+s} \la~$}~$\ze{q}{vt}f_{tm+\Hf+s}$}
  \ENDFOR
\ENDFOR
\FOR{$s = m_1, \ldots, m-1$}
  \STATE{$W_{(\ph -1)m+s} \la \ze{q}{-v}f_{-m+\Hf+s}$}
\ENDFOR
\FOR{$s = 0, \ldots, m-1$}
  \STATE{$\seq{W}[um+s]{u}{0}{\ph-1} \la \fft{\seq{W}[tm+s]{t}{0}{\ph-1}}$}
\ENDFOR
\FOR{$u=0,\ldots,\ph-1$}
  \FOR{$s=1,\ldots,m-1$}
    \STATE{$W_{um+s}\la \ze{qm}{(un+v)s}W_{um+s}$}
  \ENDFOR
  \STATE{$\seq{V}[um+\ell]{\ell}{0}{m-1} \la \fft{\seq{W}[um+s]{s}{0}{m-1}}$}
\ENDFOR
\RETURN{$\seq{V}{k}{0}{\phm-1}$}
\end{algorithmic}
\end{algorithm}
\end{minipage}
\hfill
\begin{minipage}[t]{0.5\textwidth}
\centering
\begin{algorithm}[H]
\caption{\backwardCenteredInner~\Bdesc{$v$}{$p>2$}{centered complex}~\deltadesc}
\begin{algorithmic}\hyperlabel{alg:backwardCenteredInner}
\INPUT{$\seq{F}{k}{0}{pm/2-1},\param,v$}
\STATE{$\Hf\la\floor{L/2}$}
\STATE{$\ph \la \ceil{L/(2m)}$}
\STATE{$m_0\la \phm-\Hf$}
\STATE{$m_1\la L-\Hf-(\ph-1)m$}
\FOR{$u=0,\ldots,\ph-1$}
  \STATE{$\seq{W}[um+s]{s}{0}{m-1} \la \ifft{\seq{F}[um+\ell]{\ell}{0}{m-1}}$}
  \FOR{$s=1, \ldots, m-1$}
    \STATE{$W_{um+s} \la \ze{qm}{-(un+v)s}W_{um+s}$}
  \ENDFOR
\ENDFOR
\FOR{$s=0,\ldots,m-1$}
  \STATE{$\seq{W}[tm+s]{t}{0}{\ph-1} \la \ifft{\seq{W}[um+s]{u}{0}{\ph-1}}$}
\ENDFOR
\FOR{$s = 0, \ldots, m_0-1$}
  \STATE{$V_{\Hf+s} \la W_{s}$}
\ENDFOR
\FOR{$t = 0, \ldots, \ph-1$}
  \STATE{$m_2 \la (m_1-m)\delta_{t-(\ph-1)}+m$}
  \FOR{$s = m_0\delta_t, \ldots, m_2-1$}
    \STATE{$V_{(t-\ph)m+\Hf+s} \la \ze{q}{-v(t-\ph)}W_{tm+s}$}
    \STATE{$V_{tm+\Hf+s} \la \ze{q}{-vt}W_{tm+s}$}
  \ENDFOR
\ENDFOR
\FOR{$s = m_1, \ldots, m-1$}
  \STATE{$V_{-m+\Hf+s} \la \ze{q}{v}W_{(\ph-1)m+s}$}
\ENDFOR
\RETURN{$\seq{V}{j}{0}{L-1}$}
\end{algorithmic}
\end{algorithm}
\end{minipage}
\end{multalg}

\begin{multalg}[htbp]
\begin{minipage}[t]{0.49\textwidth}
\centering
\begin{algorithm}[H]
\caption{\forwardHermitianInner~\Fdesc{$v$}{$p>2$}{Hermitian}~\deltadesc}
\begin{algorithmic}\hyperlabel{alg:forwardHermitianInner}
\INPUT{$\seq{f}{j}{0}{\Hc},\param,v$}
\STATE{$\Hc\la \ceil{L/2}$}
\STATE{$e\la \floor{m/2}+1$}
\STATE{$\ph \la \ceil{L/(2m)}$}
\STATE{$n\la q/\ph$}
\STATE{$m_0\la \min\left(\phm-\Hc+1,e\right)$}
\FOR{$s = 0, \ldots, m_0-1$}
  \STATE{$W_{s} \la f_{s}$}
\ENDFOR
\FOR{$t = 0, \ldots, \ph-1$}
  \FOR{$s = m_0\delta_t, \ldots, e-1$}
    \STATE{$W_{te+s} \la \ze{q}{v(t-\ph)}\conj{f_{(\ph-t)m-s}}+$}
    \STATE{\phantom{$W_{te+s} \la$}~$\ze{q}{vt}f_{tm+s}$}
  \ENDFOR
\ENDFOR
\FOR{$s = 0, \ldots, e-1$}
  \STATE{$\seq{W}[ue+s]{u}{0}{\ph-1} \la \fft{\seq{W}[te+s]{t}{0}{\ph-1}}$}
\ENDFOR
\FOR{$u=0,\ldots,\ph-1$}
  \FOR{$s=1,\ldots,e-1$}
    \STATE{$W_{ue+s}\la \ze{qm}{(un+v)s}W_{ue+s}$}
  \ENDFOR
  \STATE{$\seq{V}[um+\ell]{\ell}{0}{m-1} \la \crfft{\seq{W}[ue+s]{s}{0}{e-1}}$}
\ENDFOR
\RETURN{$\seq{V}{k}{0}{\phm-1}$}
\end{algorithmic}
\end{algorithm}
\end{minipage}
\hfill
\begin{minipage}[t]{0.5\textwidth}
\centering
\begin{algorithm}[H]
\caption{\backwardHermitianInner~\Bdesc{$v$}{$p>2$}{Hermitian}~\deltadesc}
\begin{algorithmic}\hyperlabel{alg:backwardHermitianInner}
\INPUT{$\seq{F}{k}{0}{pm/2-1},\param,v$}
\STATE{$\Hc\la\ceil{L/2}$}
\STATE{$e\la \floor{m/2}+1$}
\STATE{$\ph \la \ceil{L/(2m)}$}
\STATE{$n\la q/\ph$}
\STATE{$m_0\la \min\left(\phm-\Hc+1,e\right)$}
\FOR{$u=0,\ldots,\ph-1$}
  \STATE{$\seq{W}[ue+s]{s}{0}{e-1} \la \icrfft{\seq{F}[um+\ell]{\ell}{0}{m-1}}$}
  \FOR{$s=1, \ldots, e-1$}
    \STATE{$W_{ue+s} \la \ze{qm}{-(un+v)s}W_{ue+s}$}
  \ENDFOR
\ENDFOR
\FOR{$s=0,\ldots,e-1$}
  \STATE{$\seq{W}[te+s]{t}{0}{\ph-1} \la \ifft{\seq{W}[ue+s]{u}{0}{\ph-1}}$}
\ENDFOR

\FOR{$s = 1, \ldots, m_0-1$}
  \STATE{$V_{s} \la  W_{s}$}
\ENDFOR
\FOR{$t = 0, \ldots, \ph-1$}
  \STATE{$V_{tm} \la \ze{q}{-vt}W_{te}$}
  \FOR{$s = (m_0-1)\delta_t+1, \ldots, m-e$}
    \STATE{$V_{tm+s} \la \ze{q}{-vt}W_{te+s}$}
    \STATE{$V_{(\ph-t)m-s} \la \ze{q}{-v(\ph-t)}\conj{W_{te+s}}$}
  \ENDFOR
\ENDFOR
\IF{$m$ is even}
  \FOR{$t = 0, \ldots, \ph$}
    \STATE{$V_{tm+e-1} \la \ze{q}{-vt}W_{e(t+1)-1}$}
  \ENDFOR
\ENDIF
  \RETURN{$\seq{V}{j}{0}{\Hc}$}
\end{algorithmic}
\end{algorithm}
\end{minipage}
\end{multalg}

\bibliographystyle{siamplain}
\bibliography{hybrid}
\end{document}

